\documentclass[12pt]{article}
\usepackage[utf8]{inputenc}

\usepackage[top=1in, bottom=1in, left=1in, right=1in]{geometry}
\usepackage{graphicx} 
\usepackage{subfig}
\usepackage{amsthm}
\usepackage[authoryear]{natbib}
\usepackage{ulem}
\usepackage{color}
\usepackage{lineno}  
\usepackage{hyperref}
\usepackage{amsmath}
\usepackage{amssymb}
\usepackage{indentfirst} 
\usepackage{amsfonts}
\usepackage{mathrsfs}
\numberwithin{equation}{section}

\theoremstyle{plain}

\numberwithin{equation}{section}

\numberwithin{equation}{section}

\newcommand\keywords[1]{\textbf{Keywords}: #1}

\title{Stability analysis of the nonlinear pendulums under stochastic perturbations}
\author{\small{Yan Luo$^1$, Kaicheng Sheng$^2$\footnote{Corresponding author. E-mail addresses:  {\it k.sheng@sdu.edu.cn} (K.Sheng)}} \\
\small{$^1$ Research Centre for Mathematics and Interdisciplinary Sciences,}\\ \small{Shandong University, Qingdao, China }\\
\small{$^2$ School of Mathematics, Shandong University, Jinan, China }
}

\date{}

\begin{document}

%\linenumbers

\maketitle

\begin{abstract}
We consider a nonlinear pendulum whose suspension point undergoes stochastic vibrations in its plane of motion. Stochastic vibrations are constructed by stochastic differential equations with random periodic solutions. Averaging over these stochastic vibrations can be simplified with ergodicity. We give a complete description of the bifurcations of phase portraits of the averaged Hamiltonian system. The bifurcation curves of the stochastic perturbed Hamiltonian system are shown numerically. Estimations between the averaged system and the exact system are calculated. The correspondence of the averaged system to the exact system is explained through the Poincar\'e return map. Studying the averaged Hamiltonian system provided important information for the exact stochastic perturbed Hamiltonian system.

\bigskip

\noindent\keywords{Nonlinear pendulum, stochastic perturbation, random periodic solution, averaging method, phase portrait.}

\end{abstract}

\section{Introduction}

Pendulum dynamics have been foundational topics in classical mechanics. A nonlinear pendulum with a vibrating suspension point is a traditional problem of the perturbation theory, it is also called the inverted pendulum problem. The upper stable equilibrium of the nonlinear pendulum with a fast vertical vibrating suspension point was discovered by \citet{steph_1, steph_2}. \citet{kapitza1951} worked on the dynamical stability of a nonlinear pendulum with a vibrating suspension point and gave valuable insights into how periodic excitations can influence stability. The geometric aspects were discussed in \citep{levi}. Generalizations for multiple-link pendulums are presented in \citep{aches, khol}. The averaging method is used to simplify the Hamiltonian of the problem \citep{AKN}. \citet{tresch, nei} described bifurcations in the phase portraits of the averaged problem for a vibrating suspension point.

Historically, the stability of the nonlinear pendulum has been analyzed under deterministic conditions, where the external influences are either controlled or negligible. However, real-world systems frequently encounter random perturbations such as environmental noise, structural vibrations, and other unpredictable forces. The study of nonlinear pendulum dynamics with stochastic perturbations is of great interest. \citet{Freidlin, aio} studied the stability of the nonlinear pendulum with single-direction stochastic perturbations. Details about the linear and nonlinear stability of the nonlinear pendulum with stochastic perturbations on its Newton's motion equation were discussed by \citet{Shaikhet2013}. 

Averaging over stochastic problems is typically applied to dynamical systems subjected to random excitations, especially when these systems exhibit complex nonlinear behaviour \citep{XU, Fredlin, Cheban, Huang2023}. The procedure of using the ergodic theory to transfer the time-averaging of the function to its expectation in the stochastic scope is detailed in \citep{Prato}. Structures of the random dynamical system were well established by \citet{rds}. 

We consider the nonlinear pendulum with a suspension point under planar arbitrary direction stochastic vibrations. Stochastic differential equations with random periodic solutions construct the stochastic vibrations. The averaging approach of \citet{burd} is applied to the stochastic perturbed system and then the ergodic theory is used to the stochastic perturbed Hamiltonian system. As a result, a corresponding deterministic averaged Hamiltonian system is obtained. We give a complete description of the bifurcations of its phase portraits. The bifurcation curves of the stochastic perturbed system are shown numerically. The estimation between the stochastic perturbed Hamiltonian and the simplified one is studied in the probabilistic sense. Relations between the averaged system and the exact system are described through the Poincar\'e return map. 

The results of this study offer broad applicability, including the enhancement of structural resilience, the advancement of biological systems modelling, and the optimization of control systems. Studying the stability analysis of pendulums with stochastic perturbations contributes to knowledge of the behaviour of dynamical systems under uncertainty.

\section{Construction of the stochastic perturbations}

Many kinds of stochastic processes can be treated as stochastic perturbations with some small coefficients. Most of them are complicated with some further analysis. We consider the stochastic processes generated from some stochastic differential equations with ergodic measures. 

We consider a probability space $\left(\Omega,\mathscr{F},\mathbb{P}\right)$, where $\mathscr{F}$ is the Borel $\sigma$-algebra of $\Omega$ and $\mathbb{P}$ is the corresponding probability measure in $\left(\Omega,\mathscr{F}\right)$. We consider a family of mappings $\left(\vartheta_t\right)_{t\in\mathbb{R}}$ of $\left(\Omega,\mathscr{F}\right)$ into itself satisfying
the following conditions:
\begin{itemize}
    \item[a)] $\left(t,\omega \right) \to \vartheta_t\left(\omega\right)$ is measurable,
    \item[b)] $\vartheta_0={1}_\Omega$,
    \item[c)] $\vartheta_{t+s}=\vartheta_t\circ \vartheta_s$.
\end{itemize}

The combination $\left(\Omega,\mathscr{F},\mathbb{P},\left( \theta_t \right)_{t\in \mathbb{R}}\right)$ is called a $\mathbb{P}$-preserving dynamical system (or a metric dynamical system) if for a subset $A$ of $\Omega$, $\vartheta_t\mathbb{P}=\mathbb{P}\left(\vartheta^{-1}_{t}A\right)=\mathbb{P}\left(A\right)$. A metric dynamical system is called ergodic if any invariant set $A$ ($\vartheta^{-1}_tA=A$) of the dynamical system has either full measure or zero measure.

Taking $\Delta=\left\{\left(t,s\right)\in \mathbb{R}^2,s\le t \right\}$, we will consider a stochastic semi-flow $u: \Delta \times \Omega \times \mathbb{R}\to \mathbb{R}$ which  satisfies for all. $\omega\in \Omega$,
\begin{equation}
\label{ut}
u\left(t,s,\omega\right)x=u\left(t,r,\omega\right)\circ u\left(r,s,\omega\right)x, \text{ for all } s\le r\le t \text{ and } x\in \mathbb{R},
\end{equation}
and $u\left(t,s,\omega\right)x:\mathbb{R}\to \mathbb{R}$ is continuous for all $\left(t,s\right)\in \Delta$.

The two-parameter random map $u$ is called a stochastic $\tau$-periodic semi-flow if it satisfies formula (\ref{ut}) and the random periodic property: for all $\left(t,s\right)\in \Delta$, there exists a constant $\tau >0$ such that 
\begin{equation*}
\label{uttau}
u\left(t+\tau,s+\tau,\omega\right)=u\left(t,s,\theta_\tau \omega \right), \text{ for almost all } \omega \in \Omega.
\end{equation*}

\citet{prp} put forward that a random $\tau$-periodic path of  of stochastic periodic semi-flows is  a measurable process $Y:\mathbb{R}\times \Omega \to \mathbb{R}$, such that for any $\left(t,s\right)\in \Delta$,
$$
u\left(t,s,\omega\right)Y\left(s,\omega\right)=Y\left(t+s,\omega\right) \text{, } Y\left(s+\tau,\omega\right)=Y\left(s,\theta_{\tau} \omega \right) \text{ for almost all } \omega \in \Omega.
$$
Since $\vartheta_s$ is $\mathbb{P}$-preserving, then for any subset $A\subset \mathbb{R}$ and for any $\left(t,s\right)\in \Delta$, we have
$$
\mathbb{P}\left(Y\left(t,\omega\right)\in A\right)=\mathbb{P}\left(Y\left(t,\vartheta_s\omega\right)\in A\right).
$$
Noting that the Markov transition probability function is given by
\begin{equation}   \label{probability}P\left(t,s,x,A\right)=\mathbb{P}\left(\omega\in \Omega: u\left(t,s,\omega,x\right)\in A\right), \,t\ge s.
\end{equation}
For $\left(t,s\right)\in\Delta$, the adjoint operator $P^{*}(t,s)$ acting on $\mathcal{P}(\mathbb{X})$  (the space of probability measures on $(\mathbb{X}, \mathscr{B})$)  is defined by
$$
\left(P^{*}(t,s) \mu\right)(\Gamma)=\int_{\mathbb{X}} P(t,s, x, \Gamma) \mu(d x), \quad  \mu \in \mathcal{P}(\mathbb{X}), \, \Gamma \in \mathscr{B}.
$$
The measure-valued function $\rho: \mathbb{R} \rightarrow \mathcal{P}(\mathbb{X})$ is called a $\tau$-periodic measure of $\tau$-periodic Markov
transition probability \citep{rpp}, if for any $\left(t,s\right)\in \mathbb{R}$,
 \begin{equation*}
P^{*}(t,s) \rho_{s}=\rho_{t}, \, \rho_{s+\tau}=\rho_{s}. 
\end{equation*}

The random periodic paths and periodic measures can be obtained from some stochastic differential equations in $\mathbb{R}$. Let $D$ be any bounded subset in $\mathbb{R}$, and define
$$
\hat{D}=\bigcup\limits_{D\subset\mathbb{R}}D.
$$
We consider random periodic paths
$$
\xi_{i}\left(t,\omega\right)=X_t\in \hat{D}, \quad i=1,2
$$
be $\mathscr{{B}}\left(\mathbb{R}^+\right)\otimes\mathscr{B}\left(\Omega\right)$ - measurable stochastic processes generated by stochastic differential equations:
\begin{equation}
\label{sde}
dX_t=b_{i}\left(t,X_t\right)dt+\beta_{i}dW_t\quad i=1,2,
\end{equation}
where $\beta_i$ are some finite positive constant, $W_t$ is a one-dimensional standard Wiener process over a probability space $\left(\Omega,\mathscr{F}, \mathbb{P}\right)$, $b_{i}\left(t,x\right):\mathbb{R}^+\times \mathbb{R}\to \mathbb{R}$ is $\tau$-periodic:
$$
b_{i}\left(t,x\right)=b_{i}\left(t+\tau,x\right),\quad 
$$
and it is also assumed to satisfy the Lipchitz conditions:
$$
|\left(b_i\left(t,x\right)-b_i\left(t,y\right)\right)|\le \alpha_i|x-y| \text{ for all } x,y\in \mathbb{R} \text{ and for some finite } \alpha_i>0,
$$

According to \citep{rpp, prp}, the random $\tau$-periodic path $ \xi_i\left(t,\omega\right)$ ($i=1,2$) can be generated by the stochastic periodic semi-flow $u$, the Markov transition probability defined by formula (\ref{probability}) and the corresponding periodic measure $\rho $ satisfy 
$$
 \frac{1}{\tau}\int_{0}^{\tau}\left(\int_{\mathbb{R}}\left(P\left(t,s,\xi,du\right)-d\rho_s\left(u\right)\right)\right)ds= 0
$$
and they are ergodic. 

We consider the stochastic perturbations defined by the integrals of the random periodic paths $ \xi_i\left(t,\omega\right)$ ($i=1,2$) over finite time. 

%For example, if equation (\ref{sde}) satisfies the $\tau$-periodic Langevin equation:
%$$
%md\dot q_t=\left(F\left(t,q_t\right)-\gamma\dot q_t\right)dt+\sigma dW_t, \quad \sigma>0,
%$$
%where $q_t\in \mathbb{R}$ is the position, constant mass $m>0$, damping constant $\gamma>0$, time-dependent force  $F\left(t,q\right)$ is Lipchitz and satisfies $F\left(t,q\right)=F\left(t+\tau,q\right)$, $W_t$ is one-dimensional Wiener process. Then \citet{eog} derived that there are random periodic paths and periodic measures. 

\section{Hamiltonian of the system}

\begin{figure}[!htbp]
  \centering
  \includegraphics[width=0.4\textwidth]{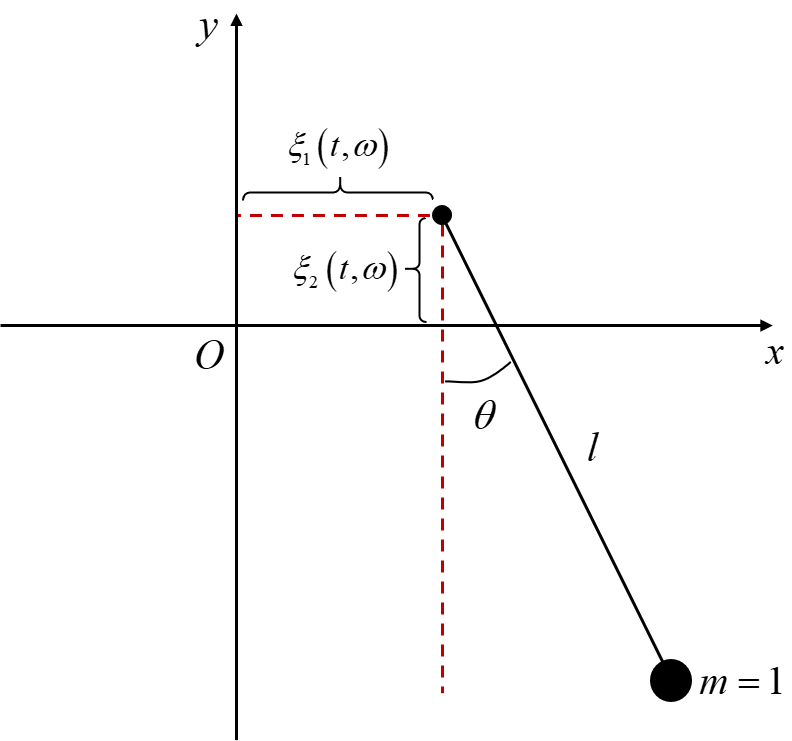}
  \caption{nonlinear pendulum under stochastic perturbation.}\label{fig1}
\end{figure}

Consider a nonlinear pendulum with a massless rod in the Cartesian frame $Oxy$, Fig. \ref{fig1}, the suspension point undergoes a pair of stochastic perturbations. Denote the angle between the pendulum rod and the vertical line straight down by $\theta$. Let $l$ and $m$ be the length of the massless rod and the mass of the bob for this pendulum. It is reasonable to set $m=1$.

We consider two random $\tau$-periodic paths $\xi_{i}\left(t,\omega\right): \mathbb{R}^+\times \Omega \to  \hat{D}$ ($i=1,2$) generated from stochastic differential equations (\ref{sde}) with zero expectations.
%corresponding to horizontal and vertical Cartesian coordinates of the pendulum's suspension point. 
%It is assumed that $\xi_{1,t}$, $\xi_{2,t}$ are independent. 

Taking small enough $\sigma_1>0$, $\sigma_2>0$,  the position of the pendulum's bob $(x, y)$ with horizontal and vertical stochastic perturbations is
\begin{equation*}
\textbf{w}=\left(x, y \right)=\left(x_0, y_0 \right)+\left(\sigma_1 \int_0^t \xi_{1}\left(s,\omega\right)\mathrm{d} s, \, \sigma_2 \int_0^t\xi_{2}\left(s,\omega\right)\mathrm{d} s \right),
\end{equation*} 
where
\begin{equation*}
\left(x_0, y_0 \right)=\left(l\cdot \sin \theta, \, -l\cdot \cos \theta \right).
\end{equation*} 
The corresponding velocity coordinate with stochastic perturbations is
\begin{equation*}
\textbf{v}=\left(\dot{x}, \dot{y} \right)=\left(l \dot{\theta } \cos \theta  ,\, l \dot{\theta }\sin \theta \right) +\left(\sigma_1 \xi_{1}\left(t,\omega\right),\, \sigma_2 \xi_{2}\left(t,\omega\right) \right).
\end{equation*} 
Then the kinetic energy of the bob is
\begin{equation*}
\begin{aligned}
    T=| \textbf{v}|=\frac{1}{2} \left( {{{\dot{x}}}^{2}}+{{{\dot{y}}}^{2}} \right)
    =& l \dot{\theta}\left(\sigma_1 \xi_{1}\left(t,\omega\right) \cos\theta+\sigma_2 \xi_{2}\left(t,\omega\right)\sin\theta \right) 
    \\
    &+\frac{1}{2} \left[l^2 {\dot\theta}^2+{\sigma_1}^2 {\xi^2_{1}\left(t,\omega\right)} +{\sigma_2}^2 {\xi^2_{2}\left(t,\omega\right)} \right],
\end{aligned}
\end{equation*} 
and the potential energy of the bob is
\begin{equation*}
V=gy=g \left[-l\,\cos\theta+\sigma_2 \cdot \int_0^t \xi_{2}\left(s,\omega\right)\mathrm{d} s \right],
\end{equation*} 
where $g$ is the gravity acceleration constant. 

Consider a pair of conjugate variables $\left(\theta,p\right)$,  the generalized momentum
\begin{equation}
\label{p}
p=\frac{\partial T}{\partial\dot \theta} = l^2 \dot\theta+l \left[\sigma_1 \xi_{1}\left(t,\omega\right) \cos\theta+\sigma_2 \xi_{2}\left(t,\omega\right)\sin\theta \right],
\end{equation} 
which is invertible at $\dot\theta$, and we also have
\begin{equation}
\label{dottheta}
\dot\theta = \frac{p}{l^2}-\frac{\sigma_1 \xi_{1}\left(t,\omega\right)\cos\theta+\sigma_2 \xi_{2}\left(t,\omega\right)\sin\theta}{l}\, .
\end{equation} 
The Lagrangian of the pendulum system is $L\left(\theta, \dot\theta, \xi_{1}, \xi_{2}\right)=T-V$. For every $t\in \mathbb{R}^+$, $\theta\in \mathbb{R}$, $\omega\in \Omega$, $L$ is a function of $\dot\theta$ which satisfies:
\begin{itemize}
    \item[a)] $L_{\left(\theta,\xi_{1},\xi_{2}\right)}\left(\dot\theta\right)\in \mathcal{C}^2 \left(\mathbb{R}\right)$,
    \item[b)] $p=\displaystyle \frac{\partial L_{\left(\theta,\xi_{1},\xi_{2}\right)}\left(\dot\theta\right)}{\partial \dot\theta}=\frac{\partial T}{\partial\dot \theta}$ is invertible at $\dot\theta$, \footnote{This condition may not be satisfied for the Lagrangian function which is not separatable of $\dot{\theta}$, if so, one may consider the corresponding stochastic differential equations and then find the existence of the solutions.}  
    \item[c)] $\displaystyle  \frac{\partial^2 L_{\left(\theta,\xi_{1},\xi_{2}\right)}\left(\dot\theta\right)}{\partial \dot\theta^2}>0$. 
\end{itemize}
%We also suppose that for all $t\in \mathbb{R}^+$, $\theta\in \mathbb{R}$, and ，
%We also suppose that for all $t\in \mathbb{R}^+, \theta\in \mathbb{R},  \omega\in \Omega$,
Therefore, for every $t\in \mathbb{R}^+$, $\theta\in \mathbb{R}$, $\omega\in \Omega$, the Legendre transformation in this  case is still applicable \citep{Arnold1974}, then the stochastic perturbed Hamiltonian for a nonlinear pendulum is  
$$
H\left( \theta,  \dot\theta, \xi_{1}, \xi_{2} \right)=p\dot{\theta}-L=p\dot{\theta}-T+V.
$$
Here $T\left( \theta,  \dot\theta, \xi_{1},\xi_2\right)={{T}_{2}}+{{T}_{1}}+{{T}_{0}}$, where ${T}_{0}$, ${T}_{1}$ and ${T}_{2}$ are terms of 1st, 2nd and 3rd order of $\dot{\theta}$. Since 
$$
p=\frac{\partial T}{\partial \dot{\theta}}=\frac{\partial {{T}_{2}}}{\partial \dot{\theta}}+\frac{\partial {{T}_{1}}}{\partial \dot{\theta}}
$$ 
and 
$$
\frac{\partial {{T}_{2}}}{\partial \dot{\theta}}\cdot\dot{\theta }=2{{T}_{2}}, \quad \frac{\partial {{T}_{1}}}{\partial \dot{\theta}}\cdot\dot{\theta}={{T}_{1}},
$$
we have
$$
p\dot{\theta }=\left( \frac{\partial {{T}_{2}}}{\partial \dot{\theta }}+\frac{\partial {{T}_{1}}}{\partial \dot{\theta }} \right)\dot{\theta }=\frac{\partial {{T}_{2}}}{\partial \dot{\theta }}\dot{\theta }+\frac{\partial {{T}_{1}}}{\partial \dot{\theta}}\dot{\theta}=2{{T}_{2}}+{{T}_{1}}.
$$
Therefore 
$$
H=2{{T}_{2}}+{{T}_{1}}-{{T}_{2}}-{{T}_{1}}-{{T}_{0}}+V={{T}_{2}}-{{T}_{0}}+V.
$$
Here $T_{0}$ is not important, ${{T}_{2}}=m{{l}^{2}}{{\dot{\theta}}^{2}}/2$. Thus it is practical to take
\begin{equation}
\label{H00}
H\left( \theta, \dot{\theta}, \xi_{1}, \xi_{2} \right)=T_{2}+V=\frac{1}{2}{{l}^{2}}{{\dot{\theta }}^{2}}+g \left[-l\,\cos\theta+\sigma_2 \cdot \int_0^t \xi_2\left( s,\omega \right)\mathrm{d} s \right].
\end{equation}
Term $ g\,\sigma_2 \cdot \int_0^t \xi_2\left( s,\omega \right)\mathrm{d} s $ does not contain $\theta$, therefore is not important here. Substitute equation (\ref{dottheta}) in to (\ref{H00}), we have 
\begin{equation}
\label{H0}
\begin{aligned}
H=\frac{1}{2}\left[\begin{aligned} &\frac{{{p}^{2}}}{{{l}^{2}}}-2\,\frac{p\left(\sigma_1 \xi_1\left( t,\omega \right) \cos\theta+\sigma_2 \xi_2\left( t ,\omega\right)\sin\theta\right)}{l} \\
&+{{\left(\sigma_1 \xi_1\left( t,\omega \right) \cos\theta+\sigma_2 \xi_2\left( t,\omega \right)\sin\theta\right)}^{2}} \end{aligned}\right] -gl\cos\theta.
\end{aligned}
\end{equation}
In fact, we observed that the Hamiltonian in the stochastic sense in this problem can still be represented by the sum of the kinetic and the potential energy of the system.  

\section{Time-averaged Hamiltonian of the system}

%In the system of Hamilton's equations  
%$$
%\dot \ph = \partial H/\partial p, \quad \dot p = - \partial H/\partial \ph \,.
%$$ 
%The right-hand side is a fast-oscillating function of time. 
In line with the averaging method \cite{bm}, for an approximate description of the dynamics of variables $(\theta, p)$, we average the stochastic perturbed Hamiltonian $H$ over time $t$. 

With the ergodicity of the corresponding probability $\mathbb{P}$, we could transfer the time-average of $\xi_i$ to the space-average which is exactly the expectation of $\xi_i$:
\begin{equation}
\label{ergodic}
\lim_{T\to \infty}\frac{1}{T}\int_{0}^{T} \xi_i\left(t,\omega\right)\mathrm{d} t =  \mathbb{E}\left(\xi_i\right)=0, \qquad i=1,\,2.
\end{equation}

The limit  $T\to \infty$ above can be omitted if we consider the system under the probability measure $\mathbb{P}$. It is because $\xi_1, \xi_2$ are random $\tau$-periodic paths and $\vartheta_t$ is $\mathbb{P}$-preserving, then for $A\subset \mathbb{R}$ and $t\in\mathbb{R}^+$, we have
\begin{equation*}
\mathbb{P}\left(\xi_i\left(t+\tau,\omega\right)\in A\right)=\mathbb{P}\left(\xi_i\left(t+\tau,\vartheta_\tau\omega\right)\in A\right)=\mathbb{P}\left(\xi_i\left(t,\omega\right)\in A\right).
\end{equation*}

The linear term of the stochastic perturbation in formula (\ref{H0}) is
$$
-\frac{p\left(\sigma_1 \xi_1\left( t ,\omega\right) \cos\theta+\sigma_2 \xi_2\left( t ,\omega\right)\sin\theta\right)}{l},
$$
which is $0$ after the averaging procedure (\ref{ergodic}). Therefore the stochastic perturbed Hamiltonian (\ref{H0}) becomes 
\begin{equation}
\label{H}
\begin{aligned}
\tilde{H}=&\frac{{{p}^{2}}}{{2 {l}^{2}}}+\frac{1}{2}{{\left[\sigma_1 \xi_1\left( t ,\omega\right) \cos\theta+\sigma_2 \xi_2\left( t ,\omega\right)\sin\theta\right]}^{2}}  -g\,l\,\cos\theta.
\\
=&\frac{{{p}^{2}}}{{2 {l}^{2}}}+\frac{1}{2} {\sigma_1^2} {\xi_1^2\left( t,\omega \right)} \cos^2\theta +{\sigma_1}{\sigma_2}{\xi_1\left( t,\omega \right)}{\xi_2\left( t ,\omega\right)} \cos\theta\sin\theta
\\
&+\frac{1}{2} {\sigma_2^2} {\xi_2^2\left( t,\omega \right)} \sin^2\theta-gl\cos\theta.
\end{aligned}
\end{equation}

For the averaging of the quadratic part of the stochastic perturbation in the Hamiltonian (\ref{H}), consideration of the averaging of ${\xi_i\left( t,\omega\right)}{\xi_j\left( t,\omega\right)}$ ($i,j=1,2$) is necessary, i.e.
$$
\lim_{T\to \infty}\frac{1}{T}\int_{0}^{T} {\xi_i\left( t ,\omega\right)}{\xi_j\left( t ,\omega\right)}\mathrm{d} t, \quad i,\,j=1,\,2.
$$

For $i=j$, taking $\eta_i\left(t,\omega\right)=\xi^2_i\left(t,\omega\right) \in \hat{D}^+$, $i=1,2$, with the ergodicity, we have
$$
\lim_{T\to \infty}\frac{1}{T}\int_{0}^{T} \eta_i\left(t,\omega\right)\mathrm{d} t =  \mathbb{E}\left(\eta_i\right)=C_{i}, \quad i=1,\,2, \, C_{i}>0.
$$

For $i\ne j$, note that 
\begin{equation*}
\begin{aligned}
&\lim_{T\to \infty}\frac{1}{T}\int_{0}^{T} {\xi_1\left( t,\omega \right)}{\xi_2\left( t ,\omega\right)}\mathrm{d} t 
=\mathbb{E}\left(\xi_1\xi_2\right)=C_{12},\quad C_{12}>0.
\end{aligned}
\end{equation*}
%We consider that 
%$$
%\zeta_{1}\left( t,\omega \right)=\left({\xi_1\left( t,\omega \right)}+{\xi_2\left( t ,\omega\right)}\right)
%$$
%and 
%$$
%\zeta_{2}\left( t,\omega \right)=\left({\xi_1\left( t,\omega \right)}-{\xi_2\left( t ,\omega\right)}\right)
%$$
%satisfy the random periodic property, then by the ergodicity of periodic measure, we have
%\begin{equation*}
%\begin{aligned}
%&\lim_{T\to \infty}\frac{1}{T}\int_{0}^{T} {\xi_1\left( t,\omega \right)}{\xi_2\left( t ,\omega\right)}\mathrm{d} t
%\\
%=&\frac{1}{2}\left(\mathbb{E}\left(\zeta_{1}\left( t,\omega \right)^2\right)-\mathbb{E}\left(\zeta_{2}\left( t,\omega \right)^2\right)\right)=C_2,\quad C_2>0.
%\end{aligned}
%\end{equation*}

Consequently, we successfully performed the averaging on the stochastic perturbed Hamiltonian (\ref{H0}), the averaged Hamiltonian is 
\begin{equation}
\label{averH}
\bar{H}=\frac{p^2}{2 l^2}+\frac{1}{2} {\sigma_1^2} C_{1} \cos^2\theta +{\sigma_1}{\sigma_2}C_{12} \cos\theta\sin\theta +\frac{1}{2} {\sigma_2^2} C_{2} \sin^2\theta-gl\cos\theta
\end{equation}
By taking 
$$
\Lambda_{1}=\frac{1}{2}(\sigma_1^2C_{1}-\sigma_2^2C_{2}) , \quad \Lambda_{2}=\frac{1}{2}{\sigma_1}{\sigma_2}C_{12}, 
$$ 
the averaged Hamiltonian up to an additive constant can be simplified to
\begin{equation*}
\label{averH1}
\bar{H}\left(\theta, p\right)=\frac{p^2}{2 l^2}+\bar{U}\left(\theta\right),
\end{equation*}
where 
\begin{equation*}
\label{barU}
\bar{U}\left(\theta\right)=\Lambda_{1} \cos 2\theta+\Lambda_{2}\sin 2\theta  - gl\cos\theta
\end{equation*}
is the effective potential of the averaged Hamiltonian. For any $\theta \in \mathbb{R}$ and any values of $\Lambda_{1}$, $\Lambda_{2}$, the function $\bar{U}\left(\theta\right)$ satisfies
\begin{itemize}
    \item[a)] $\bar{U}\left(-\theta, \Lambda_{1}, -\Lambda_{2}\right)=\bar{U}\left(\theta, \Lambda_{1}, \Lambda_{2}\right)$
    \item[b)] $\bar{U}\left(\pi-\theta, -\Lambda_{1}, \Lambda_{2}\right)=\bar{U}\left(\theta, \Lambda_{1}, \Lambda_{2}\right)$
\end{itemize}
Thus, without loss of generality, we assume that $\Lambda_{1}\ge 0$ and $\Lambda_{2} \ge 0$.

\section{Bifurcations and phase portraits}

The extrema of $\bar{U}\left(\theta\right)$ determines the stability of the nonlinear pendulum in the averaged sense. In this situation, the stochastic perturbations are approximately described by some deterministic values. Therefore, according to \citep{nei}, we partition the parameter plane $\left(\Lambda_{1}, \Lambda_{2}\right)$ of the problem into domains corresponding to different types of phase portraits of the averaged system. The bifurcation curves correspond to a) the degenerate equilibria and b) the saddle equilibria with equal values of potential energy are obtained by
\begin{itemize}
    \item[a)] $\displaystyle \frac{\partial {\bar{U}\left(\theta\right)}}{\partial \theta} = 0$, $\displaystyle \frac{\partial^2 {\bar{U}\left(\theta\right)}}{\partial \theta^2} = 0$.
    \item[b)] $\displaystyle \frac{\partial \bar{U}\left(\theta_{1}\right)}{\partial \theta_{1}}=0$, $\displaystyle \frac{\partial \bar{U}\left(\theta_{2}\right)}{\partial \theta_{2}}=0$, $\displaystyle \bar{U}\left(\theta_{1}\right)=\bar{U}\left(\theta_{2}\right)$, $\theta_{1}$ and $\theta_{2}$ are saddle equilibria.
\end{itemize}

\begin{figure}[!htbp]
\centering
\subfloat[In averaged system.]{\includegraphics[width=0.49\textwidth]{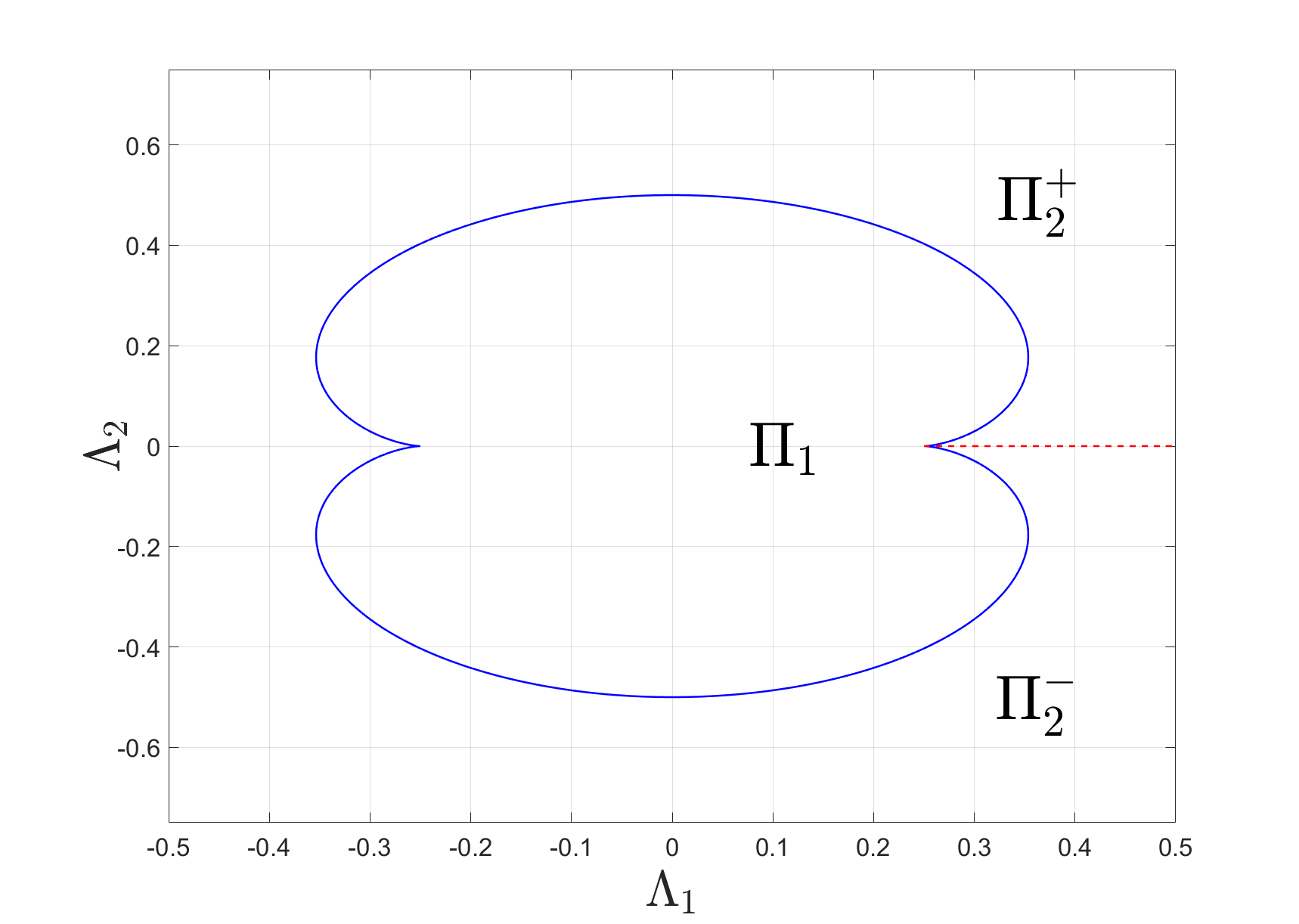}\label{Bifur1}}
\subfloat[In perturbed system.]{\includegraphics[width=0.49\textwidth]{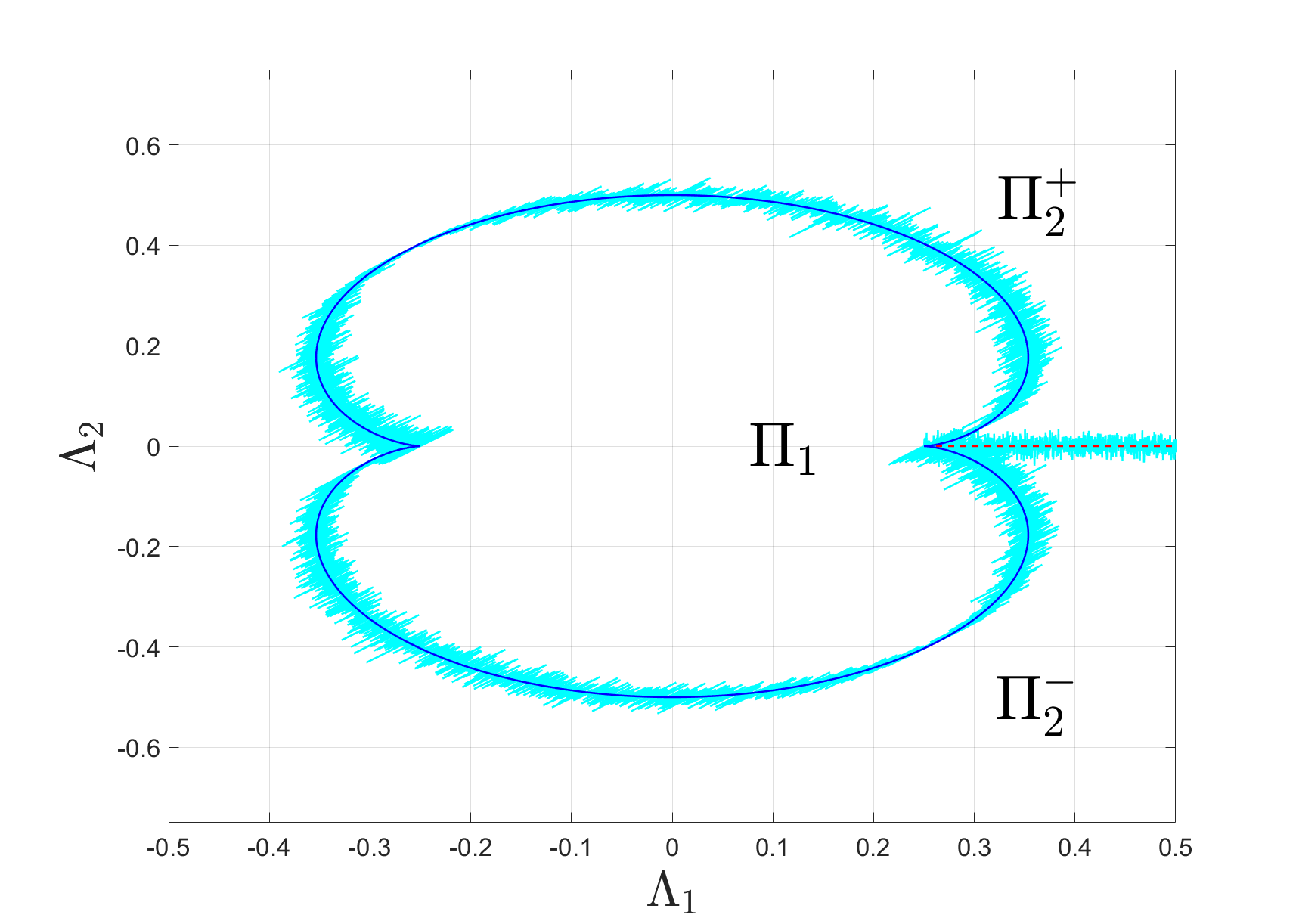}\label{Bifur2}}
\caption{Bifurcation curves}
\label{Bifur}
\end{figure}

Consequently, taking $l=1$, $g=1$, the bifurcation curves for the averaged system are the curve
\begin{equation*}
\Gamma_1 = \left\{\Lambda_{1}=\frac{\cos^3\theta}{2} -\frac{3\cos\theta}{4}, \Lambda_{2}=\frac{\sin^3\theta}{2}\right\},
\end{equation*}
and the ray 
$$
\Gamma_2 =\{\Lambda_{1}>0.25, \Lambda_{2}=0\}.
$$
The curves are shown in Fig. \ref{Bifur1}. 

\begin{figure}
\centering
\subfloat[Domain $\Pi_1$: $(\Lambda_{1}=0,\Lambda_{2}=0)$]{\includegraphics[width=0.48\textwidth]{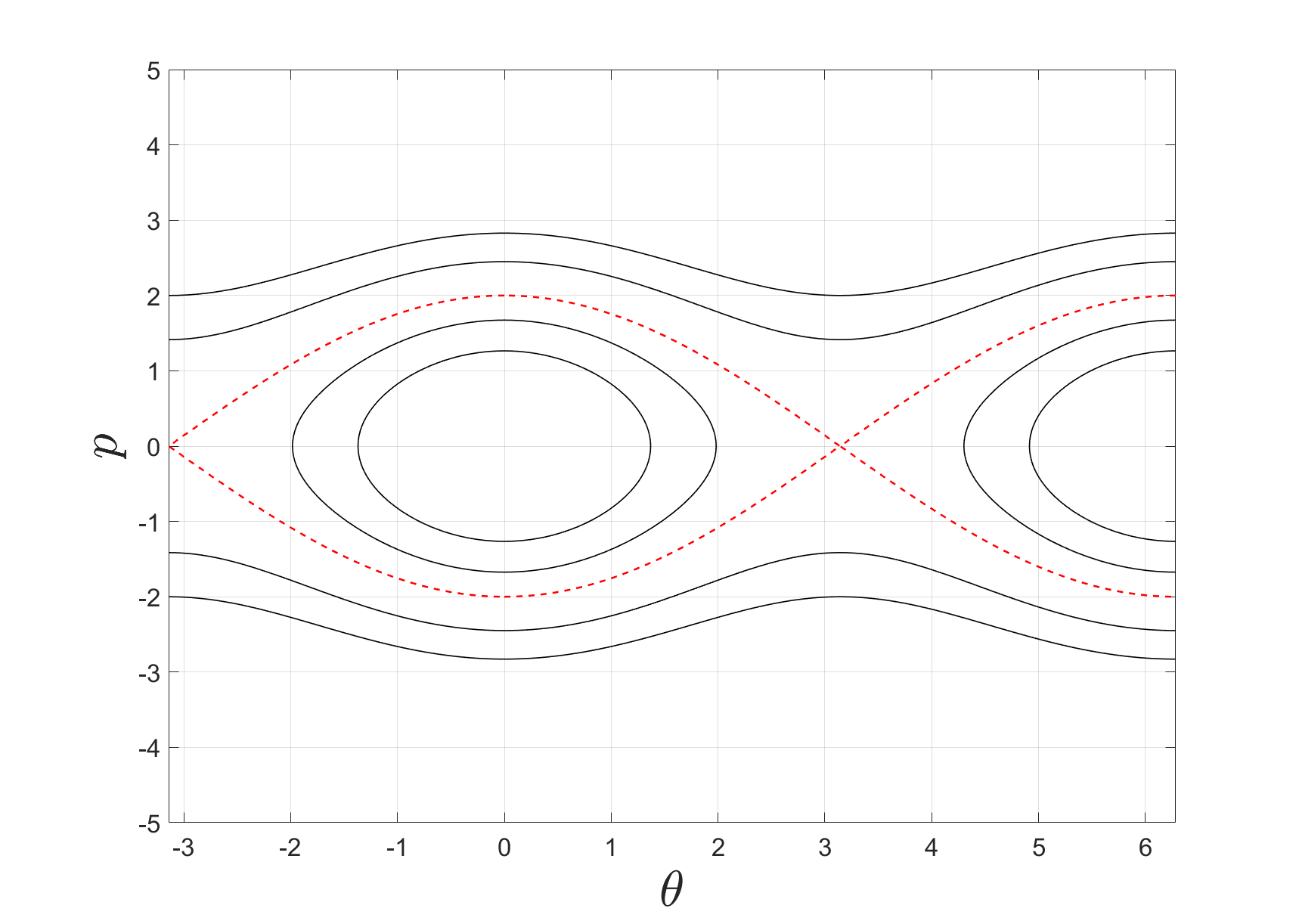}\label{phase1}}
\subfloat[Domain $\Pi_1$: $(\Lambda_{1}=0.2,\Lambda_{2}=0.1)$]{\includegraphics[width=0.48\textwidth]{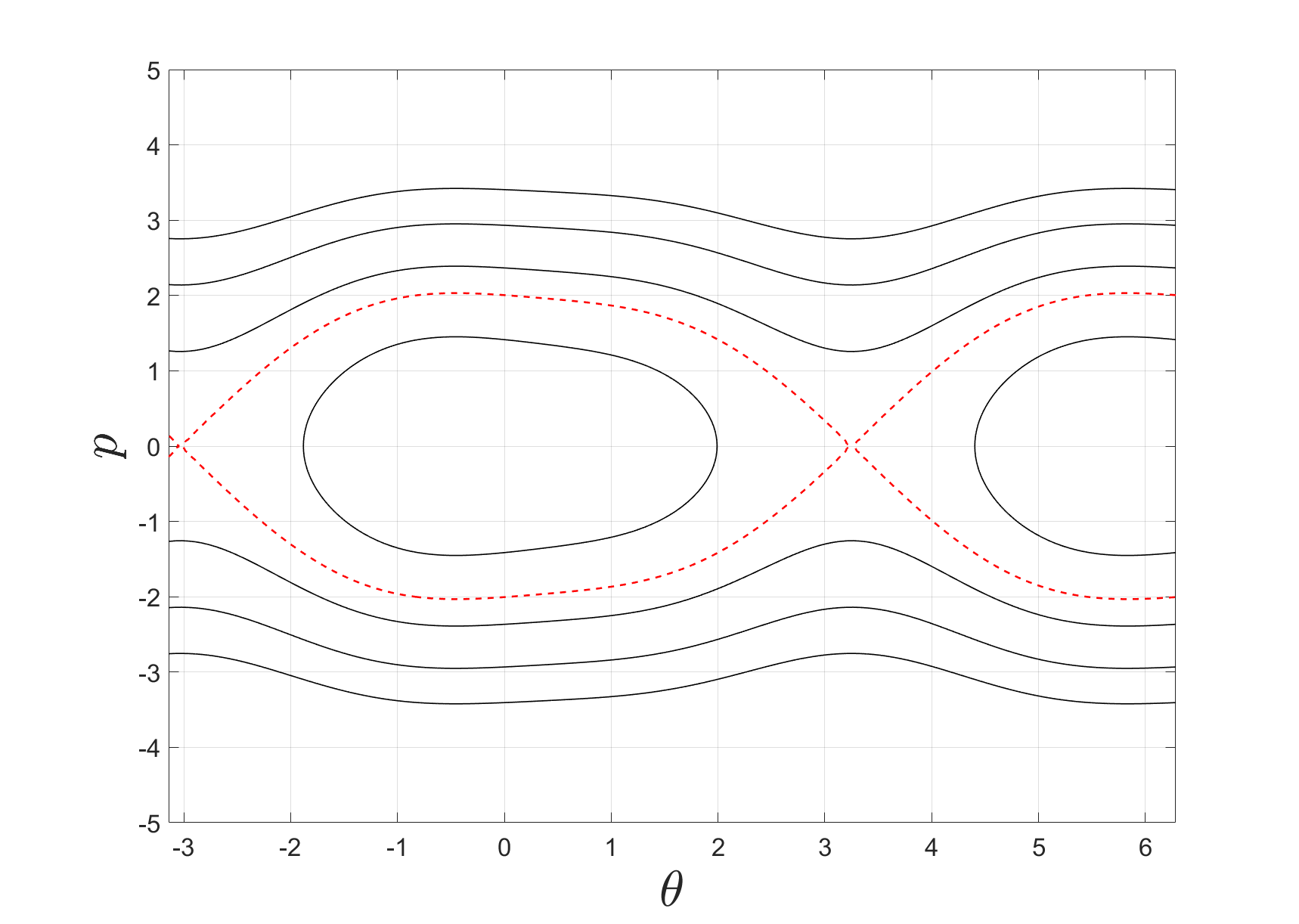}\label{phase2}}
\\
\subfloat[Domain $\Pi_2$: $(\Lambda_{1}=0.5,\Lambda_{2}=1)$]{\includegraphics[width=0.48\textwidth]{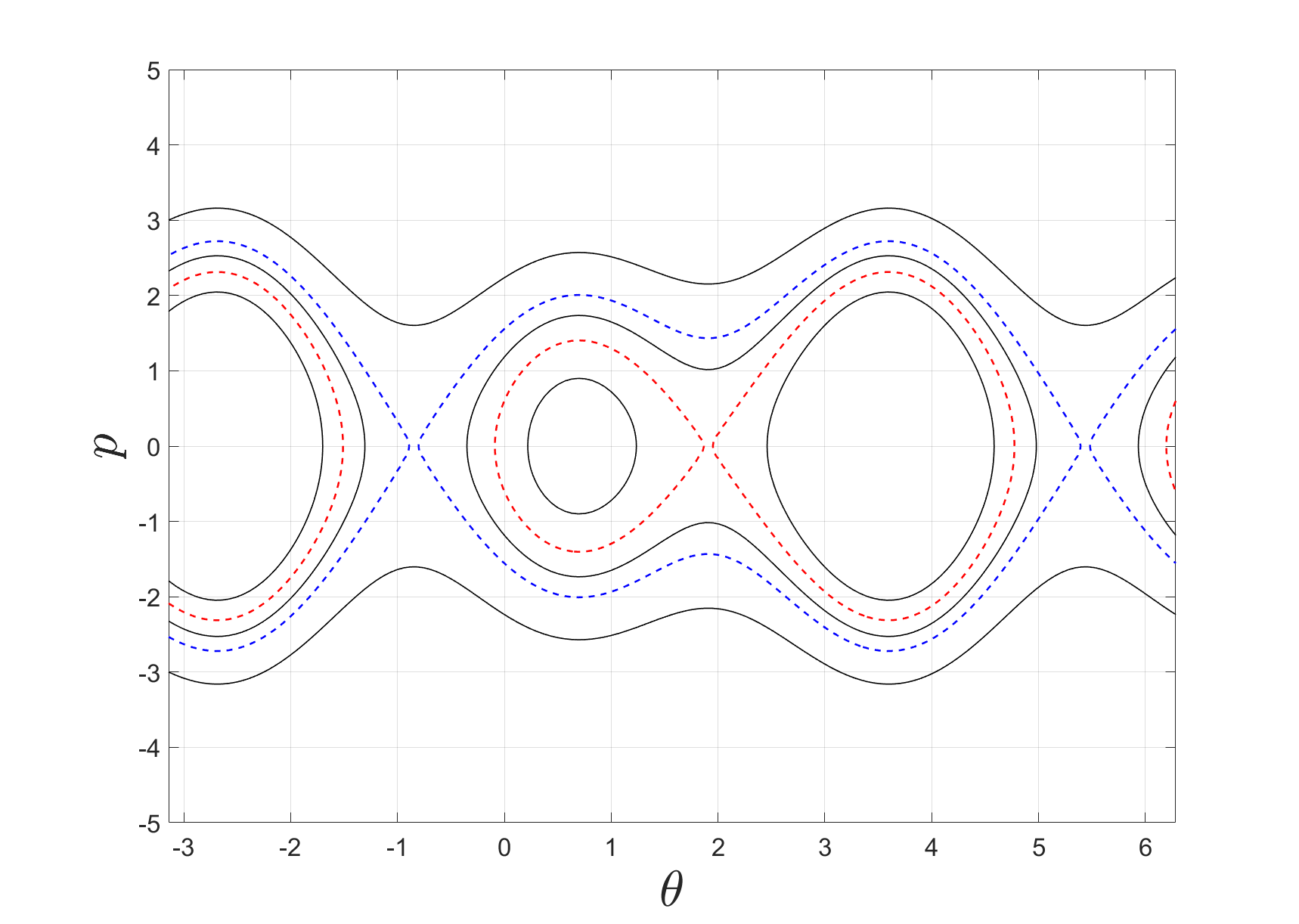}\label{phase3}}
\subfloat[Domain $\Pi_2$: $(\Lambda_{1}=0.5,\Lambda_{2}=-1)$]{\includegraphics[width=0.48\textwidth]{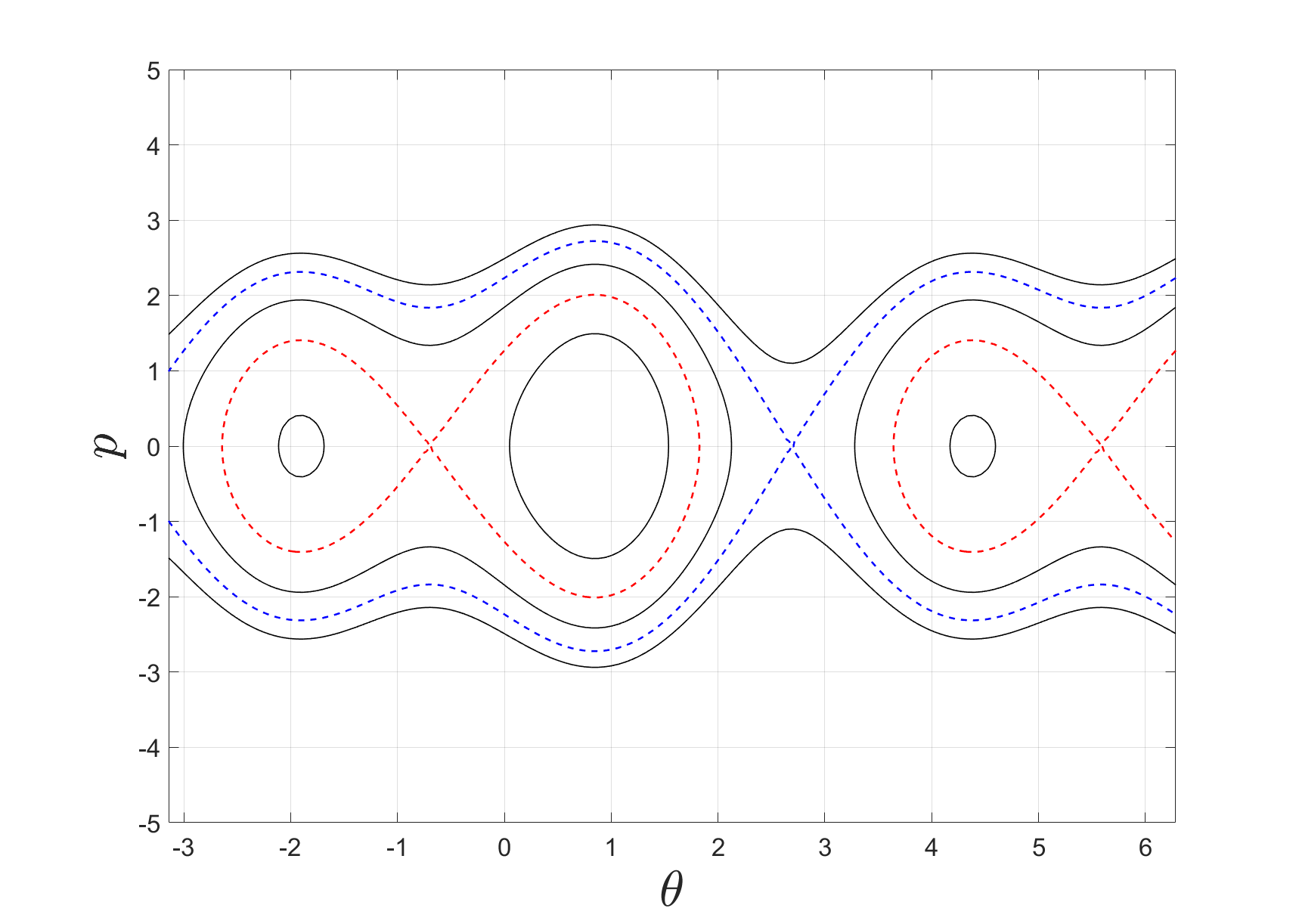}\label{phase4}}
\caption{Phase portraits}
\label{phase}
\end{figure}

In fact, due to the stochastic perturbations, the exact dynamics near the bifurcation curves are complicated. We consider the time projection trajectories of $(\Lambda_{1},\Lambda_{2})$ in the plane $O\Lambda_{1}\Lambda_{2}$, see Fig. \ref{Bifur2}, random shifts are taken place on the bifurcation curves. To estimate the shifts, we consider the effective potential in equation (\ref{H}), denoted by
\begin{equation*}
\label{tildeU}
\begin{aligned}
\tilde{U}\left( t, \theta \right)=&\frac{1}{2} {\sigma_1^2} {\xi_1^2\left( t,\omega \right)} \cos^2\theta +{\sigma_1}{\sigma_2}{\xi_1\left( t,\omega \right)}{\xi_2\left( t ,\omega\right)} \cos\theta\sin\theta
\\
&+\frac{1}{2} {\sigma_2^2} {\xi_2^2\left( t,\omega \right)} \sin^2\theta-gl\cos\theta.
\\
=&\frac{1}{2}\left({\sigma_1^2} {\xi_1^2}-{\sigma_2^2} {\xi_2^2}  \right)\cos 2\theta+ \frac{1}{2}{\sigma_1}{\sigma_2}{\xi_1}{\xi_2}\sin 2\theta  - gl\cos\theta
\end{aligned}
\end{equation*}
Then by the averaging principle, 
\begin{equation*}
%\label{es1}
\mathbb{E}\left(\left|\bar{U}\left( \theta \right)-\tilde{U}\left( t, \theta \right)\right|\right) < O(\max \{\sigma_1, \sigma_2\}).
\end{equation*}
Moreover, it is easy to check 
\begin{equation*}
%\label{es2}
\mathbb{E}\left(\frac{\partial }{\partial \theta}\left|\bar{U}\left( \theta \right)-\tilde{U}\left( t, \theta \right)\right|\right)  < O(\max \{\sigma_1, \sigma_2\}),
\end{equation*}
and 
\begin{equation*}
%\label{es3}
\mathbb{E}\left(\frac{\partial^2 }{\partial \theta^2}\left|\bar{U}\left( \theta \right)-\tilde{U}\left( t, \theta \right)\right| \right) < O(\max \{\sigma_1, \sigma_2\}).
\end{equation*}
%Detailed principle of the estimations (\ref{es1}), (\ref{es2}) and (\ref{es3}) will be shown in the next section. 
This reveals that bifurcation curves of equation (\ref{H}) in the plane $O\Lambda_{1}\Lambda_{2}$ are close to bifurcation curves $\Gamma_1$ and $\Gamma_2$ with small random shifts. Choosing different representative points in domain $\Pi_1$ and $\Pi_2$ of Fig. \ref{Bifur} respectively, the phase portraits are shown in Fig. \ref{phase}. 

It is shown that there are 2 equilibria when $\Lambda_{1}$ and $\Lambda_{2}$ are taken in Domain $\Pi_1$, one of them is stable while the other one is unstable (Fig. \ref{phase1}, \ref{phase2}). There are 4 equilibria when $\Lambda_{1}$ and $\Lambda_{2}$ are taken in Domain $\Pi_2$, two of them are stable, while the other two are unstable (Fig. \ref{phase3}, \ref{phase4}). Fig. \ref{phase1} corresponds to the classical nonlinear pendulum without perturbation.

\section{Relations to the exact problem}

\subsection{Estimation of the Hamiltonian}

Averaging over fast phases can simplify the study of slow phases. It is an approximate calculation, the estimation of the variance is important. We take averaging on the time-depended perturbation $\xi_{i}\left(t,\omega\right)$ ($i=1,2$) and its derivatives, the stochastic perturbed Hamiltonian for a nonlinear pendulum turns to be a deterministic one, the deviation between them are small for the long-time behaviour. 

The stochastic process $\xi_i$  generated by stochastic differential equation (\ref{sde}) can be considered as
\begin{equation*}
 \xi_i\left(t,\omega\right)=z_i+\int_0^t b_i\left(s,\xi_i\left(s,\omega\right)\right)ds+  \int_0^t \beta_idW_s.\quad i=1,2, 
\end{equation*}
where $z_i<\infty$ is the initial value of $\xi_i$. And from formulas (\ref{p}),(\ref{H0})  and (\ref{averH}), we have
\begin{equation*}
%\label{H-barH}
    \begin{aligned}
        &|H-\bar{H}|
        \\
        =&\left|l\dot\theta+\left(\sigma_1\xi_1\cos \theta+\sigma_2\xi_2\sin \theta\right)^2+\frac{1}{2}\left(\sigma^2_1C_1\cos^2\theta+2\sigma_1\sigma_2 C_{12}\cos\theta\sin\theta+\sigma^2_2 C_2\sin^2\theta\right)\right|
        \\
        \le&|l\dot\theta|+|\sigma_1\xi_1|^2+2|\sigma_1\sigma_2\xi_1\xi_2|+|\sigma_2\xi_2|^2+\frac{1}{2}|\sigma^2_1C_1|+2|\sigma_1\sigma_2 C_{12}|+|\sigma^2_2 C_2|.
    \end{aligned}
\end{equation*}
Taking 
\begin{equation*}
\begin{aligned}
M_1&=|\sigma_1\xi_1|^2+2|\sigma_1\sigma_2\xi_1\xi_2|+|\sigma_2\xi_2|^2, 
\\
M_2&=|\frac{1}{2}|\sigma^2_1C_1|+2|\sigma_1\sigma_2 C_{12}|+|\sigma^2_2 C_2|,
\end{aligned}
\end{equation*}
then for any $\delta>0$, 
\begin{equation}
\label{H-barH}
\mathbb{P}\left(|H-\bar{H}|>\delta\right)\le \mathbb{P}\left(|l\dot\theta|+M_1 +M_2 >\delta \right).
\end{equation}
Here,
\begin{equation*}
\mathbb{P}\left(|l\dot\theta|+M_1 +M_2 >\delta \right)=\mathbb{P}\left(M_1>\delta-|l\dot\theta|-M_2\right)=\mathbb{P}\left(M_1>\hat{\delta}\right), 
\end{equation*}
where $\hat{\delta}=\delta-|l\dot\theta|-M_2$. 
By Chebyshev’s inequality, we have for any $\hat{\delta}>0$,
\begin{equation*}
\begin{aligned}
\mathbb{P}\left(M_1>\hat{\delta}\right)\le&\frac{1}{\hat{\delta}^2} \mathbb{E}\left(M^2_1\right)\\
\le &\frac{1}{\hat{\delta}^2}\left(\mathbb{E}\left(|\sigma_1\xi_1|^4\right)+\mathbb{E}\left(|\sigma_2\xi_2|^4\right)+4\mathbb{E}\left(|\sigma_1\sigma_2\xi_1\xi_2|^2\right)\right)\\
&+\frac{2}{\hat{\delta}^2}\left(\mathbb{E}\left(2|\sigma^3_1\sigma_2\xi^3_1\xi_2|\right)+\mathbb{E}\left(|\sigma^2_1\sigma^2_2\xi^2_1\xi^2_2|\right)+\mathbb{E}\left(2|\sigma_1\sigma^3_2\xi_1\xi^3_2\right)\right)
 \end{aligned}
\end{equation*}

Now by the Lipchitz conditions in Section 2 and Itô formula, for a given $t\in [0,\tau]$, $i=1,2$, we have
\begin{equation*}
    \label{xi12}
    \begin{aligned}  
    |\xi^2_i\left(t,\omega\right)|\le&|z_i|^2+ 2\int_0^t |\xi_i b_i\left(s,\xi_i\right)|ds +\int_0^t |\beta^2_i|ds +\int_0^t2|\xi_i \beta_i|dW_s
        \\   
        \le&|z_i|^2+\int_0^t(2\alpha_i\xi^2_i+\beta^2_i)ds+\int_0^t2\beta_i|\xi_i|dW_s.
\\
    |\xi^3_i \left(t,\omega\right)|\le & |z_i|^3+3\left(\int_0^t\left(|\xi^2_ib_i\left(s,\xi_i\right)|+|\xi_i\beta^2_i|\right)ds+\int_0^t|\xi^2_i\beta_i|dW_s\right)
    \\    
    \le & |z_i|^3+3\left(\int_0^t\left(\alpha_i|\xi^3_i|+\beta^2_i|\xi_i|\right)ds+\int_0^t\beta_i|\xi^2_i|dW_s\right). 
    \end{aligned}
\end{equation*}
Then for $i=1,2$, we have 
\begin{equation}
\label{xi2}
    \begin{aligned}
&\mathbb{E}\left(|\sigma_i\xi_i|^4\right)\\
\le&\sigma^4_i\mathbb{E}\left(\int_0^t \left(2\alpha_i\xi^2_i\left(s,\omega\right)+\beta^2_i\right)ds\right)^2+\sigma^4_i\mathbb{E}\left(\int_0^t2\beta_i|\xi_i\left(s,\omega\right)|dW_s\right)^2\\
&+\sigma^4_i|z_i|^4+2\sigma^4_i|z_i|^2\mathbb{E}\left(\int_0^t \left(2\alpha_i\xi^2_i\left(s,\omega\right)+\beta^2_i\right)ds\right)\\
\le&\sigma^4_i\mathbb{E}\left(\int_0^t \left(2\alpha_i\xi^2_i\left(s,\omega\right)+\beta^2_i\right)^2ds\right)+\sigma^4_i\mathbb{E}\left(\int_0^t4\beta^2_i\xi^2_i\left(s,\omega\right)ds\right)\\
&+\sigma^4_i|z_i|^4+2\sigma^4_i|z_i|^2\mathbb{E}\left(\int_0^t \left(2\alpha_i\xi^2_i\left(s,\omega\right)+\beta^2_i\right)ds\right)\\
\le&\sigma^4_i\mathbb{E}\left(\int_0^t \left(4\alpha^2_i\xi^4_i\left(s,\omega\right)+\beta^4_i+4\alpha_i\beta^2_i\xi^2_i\left(s,\omega\right)+4\beta^2_i\xi^2_i\left(s,\omega\right)\right)ds\right)\\
&+\sigma^4_i|z_i|^4+2\sigma^4_i|z_i|^2\mathbb{E}\left(\int_0^t \left(2\alpha_i\xi^2_i\left(s,\omega\right)+\beta^2_i\right)ds\right)\\
\le&\sigma^4_it\hat{C}_i+\sigma^4_i|z_i|^4,
    \end{aligned}
\end{equation}
where $\hat{C_i}>0$ ($i=1,2$) are some finite constants. Similarly, 
\begin{equation}
\label{xi1xi2}
    \begin{aligned}
&\mathbb{E}\left(\sigma^2_1\sigma^2_2\xi^2_1\xi^2_2\right)\\ 
        %&\le \sigma^2_1\sigma^2_2\mathbb{E}\left(\left(\int_0^t\left(2\alpha_1+\beta^2_1\right)\xi^2_1\left(s,\omega\right)ds+\int_0^t2\beta_1\xi^2_1dW_s\right)\left(\int_0^t\left(2\alpha_2+\beta^2_2\right)\xi^2_2\left(s,\omega\right)ds+\int_0^t2\beta_2\xi^2_2dW_s\right)\right)\\
        \le &\sigma^2_1\sigma^2_2 \mathbb{E}\left(\int_0^t\left(4\alpha_1\alpha_2\xi^2_1\xi^2_2+2\alpha_1\beta^2_2\xi^2_1+2\alpha_2\beta^2_1\xi^2_2+\beta^2_1\beta^2_2\right)ds+\int_0^t4\beta_1\beta_2|\xi_1\xi_2|ds\right)\\
&+\sigma^2_1\sigma^2_2|z_1z_2|+\sigma^2_1\sigma^2_2|z_1|\mathbb{E}\left(\int_0^t\left(2\alpha_2\xi^2_2\left(s,\omega\right)+\beta^2_2\right)ds\right)\\
&+\sigma^2_1\sigma^2_2|z_2|\mathbb{E}\left(\int_0^t\left(2\alpha_1\xi^2_1\left(s,\omega\right)+\beta^2_1\right)ds\right)\\
\le&\sigma^2_1\sigma^2_2t\hat{C}_{12}+\sigma^2_1\sigma^2_2|z_1z_2|,
        %&\frac{16\sigma^2_1\sigma^2_2\left(4\alpha_1\alpha_2+2\alpha_1\beta^2_2+2\alpha_2\beta^2_1+\beta^2_1\beta^2_2+4\beta_1\beta_2\right)\tau}{\delta^2}\hat{C}_{12}
    \end{aligned}
\end{equation}
%Define $K=:max\{l\dot\theta\}$, $\hat{K}=\frac{1}{2}|\sigma^2_1C_1+2\sigma_1\sigma_2C_2+\sigma^2_2C_3|$, then by formulas (\ref{xi2}) and (\ref{xi1xi2}),  
where $\hat{C}_{12}>0$ is some finite constant. Also, we have
\begin{equation}
\label{xi3}
    \begin{aligned}
&\mathbb{E}\left(|\sigma^3_i\sigma_j\xi^3_i\xi_j|\right)\\
%& \le 3\sigma^3_i\sigma_j\mathbb{E}\left(\left(\int_0^t\left(\alpha_i+\beta^2_i\right)|\xi^3_i\left(s,\omega\right)|ds+\int_0^t\beta_i|\xi^3_i\left(s,\omega\right)|dW_s\right)\left(\int_0^t\alpha_j|\xi_j|ds+\int_0^t\beta_j|\xi_j|dW_s\right)\right)\\
\le&3\sigma^3_i\sigma_j\mathbb{E}\left(\int_0^t\left(\alpha_i\alpha_j|\xi^3_i\xi_j|+\beta^2_i\alpha_j|\xi_i\xi_j|\right)ds+\int_0^t\beta_i\beta_j|\xi^2_i\xi_j|ds\right)+3\sigma^3_i\sigma_j|z^3_iz_j|\\
&+3\sigma^3_i\sigma_j|z_j|\mathbb{E}\left(\int_0^t\left(\alpha_i|\xi^3_i|+\beta^2_i|\xi_i|\right)ds\right)+3\sigma^3_i\sigma_j|z_i|^3\mathbb{E}\left(\int_0^t\alpha_j|\xi_j|ds\right)\\
\le& 3\sigma^3_i\sigma_jt\tilde{C}_{ij}+3\sigma^3_i\sigma_j|z^3_iz_j|, \text{ for some } \tilde{C}_{ij}>0.
    \end{aligned}
\end{equation}
Combining inequalities (\ref{xi2}), (\ref{xi1xi2}) and (\ref{xi3}), for $t\in [0,\tau]$, we have
$$
\lim\limits_{\sigma_1,\sigma_2\to 0}\mathbb{P}\left(M_1>\hat{\delta}\right)=0. 
$$
Therefore by (\ref{H-barH}), for any $\delta>0$, we have
\begin{equation*}
\lim\limits_{\sigma_1,\sigma_2\to 0}\mathbb{P}\left(|H-\bar{H}|>\delta\right)=0.
\end{equation*}

Consequently, we proved that on the same energy level, the manifold of the stochastic perturbed system is close to that of the averaged system. The maximum distance of two manifolds turns to zero when the largest amplitude of the stochastic perturbations goes to zero. 

\subsection{On the Poincar\'e return map}

The trajectories on the Poincar\'e return map greatly present the dynamics of the exact system. Finding the relations between these trajectories and the phase curve in the averaged system reveals the correspondence between the averaged and exact system. We will discuss the Poincar\'e return map in the random regime.

We consider the Poincar\'e section for transition probability. \citet{rpp} defined that the collection of subsets $\left\{L_s: s \geqslant 0\right\} \subset \mathcal{B}\left(\mathbb{R}^d\right)$ are called the Poincar\'e sections of the transition probability function $P(t, s,x, \cdot)$ if
$$
L_{s+\tau}=L_s,
$$
and for any $s \in[0, \tau), t \geqslant 0$,
$$
P\left(t+s, s, x, L_{s+t}\right)=1, \quad x \in L_s .
$$
With the periodic property of $L_s$, we have for a given $t\in \mathbb{R}^+$, 
$$
P\left(t+n\tau, t, x, L_{t}\right)=1, \quad x \in L_t .
$$

This means that starting from $x \in L_t, u\left(t+n\tau,t,\omega, x\right)$ returns to the set $L_t$ with probability one at any time being a multiple integral of the period. This could be regarded as the Poincaré returning map property in the random regime, mirroring the celebrated Poincaré mapping in the deterministic case. %However, the map $u\left(t+n\tau,t, \omega,x\right)$ does not have a fixed point on $L_t$. This is very different from the deterministic case.

According to the \cite{prp}, the random periodic path $ \xi_i\left(t,\omega\right)$ ($i=1,2$) can be generated by the stochastic periodic semi-flow $u$. Consider the Poincar\'e return map for plane $\theta, p(\xi_1,\xi_2)$ in the stochastic dynamical system, the Poincar\'e sections are $t=0 \ {\rm mod}\, \tau$.  Then we have: 
\begin{itemize}
    \item[a)] The random periodic trajectories of the Poincar\'e return map in the stochastic regime fill the plane $\theta, p(\xi_1,\xi_2)$ up to a small remainder corresponding to the amplitudes of the stochastic perturbations. 
    \item[b)] Supposed the equilibria in the averaged system is $E_{0}$. For a given $t\in [0,\tau]$,  $a_t\left(\xi_1,\xi_2\right)=\left(\theta,p\left(\xi_1,\xi_2\right)\right)$, there is a set $\hat{L}_t$  such that
            $$
            \mathbb{P}\left(a_t\left(\xi_1,\xi_2\right)\in \hat{L}_t\right)=1. 
            $$
    Then, for any $\delta>0$, we have 
    $$
    \lim\limits_{\sigma_1,\sigma_2\to 0}\mathbb{P}\left(|a_t\left(\xi_1,\xi_2\right)-E_0|>\delta\right)=0.
    $$
    \item[c)] There is an exponentially small splitting phenomenon between the stable and unstable trajectories of the adjacent hyperbolic sets $\hat{L}_t$, $\hat{L}_{t+\tau}$ of the map. One may use the Melnikov function to estimate this splitting, which is an open question. These trajectories are close to the separatrices of the averaged system. 
\end{itemize}

Thus, the averaged system is bound up with the exact (not averaged) problem in the dynamics. Studying the averaged Hamiltonian system provided important information for the exact stochastic perturbed Hamiltonian system.

\section*{Acknowledgement}

The authors express their gratitude to Dr. Baoyou Qu for the discussions. Dr. Yan Luo thanks the National Natural Science Foundation of China (NSFC) for the support of this research (Grant: 12471142). 

\section*{Conflict of interest statement}

The authors declare that there are no conflicts of interest, we do not have any possible conflicts of interest.

\bibliographystyle{plainnat}  
\bibliography{reference}

\begin{thebibliography}{23}
\providecommand{\natexlab}[1]{#1}
\providecommand{\url}[1]{\texttt{#1}}
\expandafter\ifx\csname urlstyle\endcsname\relax
  \providecommand{\doi}[1]{doi: #1}\else
  \providecommand{\doi}{doi: \begingroup \urlstyle{rm}\Url}\fi

\bibitem[Acheson(1993)]{aches}
D.~Acheson.
\newblock A pendulum theorem.
\newblock \emph{Proc. R. Soc. London}, 433:\penalty0 239--245, 1993.

\bibitem[Arnold(2003)]{rds}
L.~Arnold.
\newblock \emph{Random Dynamical Systems}.
\newblock Springer-Verlag, Berlin Heidelberg New York, 2003.

\bibitem[Arnold(1974)]{Arnold1974}
V.~I. Arnold.
\newblock \emph{Mathematical Methods of Classical Mechanics}.
\newblock Springer, 1974.

\bibitem[Arnold et~al.(2006)Arnold, Kozlov, and Neishtadt]{AKN}
V.~I. Arnold, V.~V. Kozlov, and A.~I. Neishtadt.
\newblock \emph{Mathematical aspects of classical and celestial mechanics}.
\newblock Springer, Berlin, 2006.

\bibitem[Bogoliubov and Mitropolsky(1961)]{bm}
N.~N. Bogoliubov and Y.~A. Mitropolsky.
\newblock \emph{Asymptotic methods in the theory of non-linear oscillations}.
\newblock Gordon and Breach Sci. Publ., New York, 1961.

\bibitem[Burd and Matveev(1985)]{burd}
V.~Sh. Burd and V.~N. Matveev.
\newblock \emph{Asymptotic methods on an infinite interval in problems of nonlinear mechanics}.
\newblock Yaroslavl' Univ., Yaroslavl', 1985.

\bibitem[Cheban and Liu(2021)]{Cheban}
D~Cheban and Z~Liu.
\newblock Averaging principle on infinite intervals for stochastic ordinary differential equations.
\newblock \emph{Electronic Research Archive}, 29\penalty0 (4):\penalty0 2791--2817, 2021.

\bibitem[{Da Prato} and Zabczyk(1996)]{Prato}
G.~{Da Prato} and J.~Zabczyk.
\newblock \emph{Ergodicity for Infinite Dimensional Systems}.
\newblock Cambridge University Press, 1996.
\newblock ISBN 9780521579001.
\newblock This book is devoted to the asymptotic properties of solutions of stochastic evolution equations in infinite dimensional spaces.

\bibitem[Feng and Zhao(2020)]{rpp}
C.~Feng and H.~Zhao.
\newblock Random periodic processes, periodic measures and ergodicity.
\newblock \emph{Journal of Differential Equations}, 269:\penalty0 7382--7428, 2020.

\bibitem[Feng et~al.(2011)Feng, Zhao, and Zhou]{prp}
C.~Feng, H.~Zhao, and B.~Zhou.
\newblock Pathwise random periodic solutions of stochastic differential equations.
\newblock \emph{Journal of Differential Equations}, 251:\penalty0 119--149, 2011.

\bibitem[Freidlin and Weber(1999)]{Freidlin}
M.~Freidlin and M.~Weber.
\newblock A remark on random perturbations of the nonlinear pendulum.
\newblock \emph{Annals of Applied Probability}, 9, 08 1999.
\newblock \doi{10.1214/aoap/1029962806}.

\bibitem[Freidlin and Wentzell(2012)]{Fredlin}
M.~Freidlin and A.~Wentzell.
\newblock \emph{Random Perturbations of Dynamical Systems}.
\newblock Springer Heidelberg New York Dordrecht London, 2012.
\newblock ISBN 9783642258466.

\bibitem[Huang et~al.(2023)Huang, Kuksin, and Piatnitski]{Huang2023}
G.~Huang, S.~Kuksin, and A.~Piatnitski.
\newblock Averaging for stochastic perturbations of integrable systems.
\newblock \emph{Preprint}, 2023.

\bibitem[Kapitza(1951)]{kapitza1951}
P.~L. Kapitza.
\newblock Dynamic stability of a pendulum when its point of suspension vibrates.
\newblock \emph{Soviet Physics Journal of Experimental and Theoretical Physics}, 21\penalty0 (3):\penalty0 588--597, 1951.

\bibitem[Kholostova(2009)]{khol}
O.~Kholostova.
\newblock On the motions of a double pendulum with vibrating suspension point.
\newblock \emph{Mechanics of Solids}, 44:\penalty0 184--197, 2009.

\bibitem[Levi(1998)]{levi}
M.~Levi.
\newblock Geometry of kapitsa's potential.
\newblock \emph{Nonlinearity}, 11:\penalty0 1365--1368, 1998.

\bibitem[Neishtadt and Sheng(2017)]{nei}
A.~I. Neishtadt and K.~Sheng.
\newblock Bifurcations of phase portraits of pendulum with vibrating suspension point.
\newblock \emph{Communications in Nonlinear Science and Numerical Simulation}, 47:\penalty0 71--80, 2017.

\bibitem[Ovseyevich(2006)]{aio}
A.~I. Ovseyevich.
\newblock The stability of an inverted pendulum when there are rapid random oscillations of the suspension point.
\newblock \emph{J. Appl. Math. Mech.}, 70:\penalty0 761--768, 2006.

\bibitem[Shaikhet(2013)]{Shaikhet2013}
L.~Shaikhet.
\newblock \emph{Lyapunov Functionals and Stability of Stochastic Functional Differential Equations}.
\newblock Springer Cham Heidelberg, 2013.
\newblock ISBN 9783319001005.

\bibitem[Stephenson(1908{\natexlab{a}})]{steph_1}
A.~Stephenson.
\newblock On induced stability.
\newblock \emph{Philosophical Magazine Series 6}, 15:\penalty0 233--236, 1908{\natexlab{a}}.

\bibitem[Stephenson(1908{\natexlab{b}})]{steph_2}
A.~Stephenson.
\newblock On a new type of dynamical stability.
\newblock \emph{Memoirs and Proceedings of the Manchester Literary and Philosophical Society}, 52\penalty0 (8):\penalty0 1--10, 1908{\natexlab{b}}.

\bibitem[Treschev(2008)]{tresch}
D.~Treschev.
\newblock Some aspects of finite-dimensional hamiltonian dynamics.
\newblock In W.~Craig, editor, \emph{Hamiltonian dynamical systems and applications}, pages 1--19. Springer, Dordrecht, 2008.

\bibitem[Xu et~al.(2011)Xu, Duan, and Xu]{XU}
Y.~Xu, J~Duan, and W.~Xu.
\newblock An averaging principle for stochastic dynamical systems with lévy noise.
\newblock \emph{Physica D: Nonlinear Phenomena}, 240\penalty0 (17):\penalty0 1395--1401, 2011.

\end{thebibliography}

\end{document}